\documentclass{elsarticle}

\usepackage{hyperref,microtype,amssymb,mathtools}
\usepackage{algorithm}
\usepackage[noend]{algpseudocode}

\usepackage[nolist]{acronym}
\newcommand{\ie}{\emph{i.e.}}
\newcommand{\eg}{\emph{e.g.}}
\newcommand{\requests}{R}
\newcommand{\locations}{L}

\newcommand{\origin}[1]{#1^o}
\newcommand{\destination}[1]{#1^d}
\newcommand{\tworigin}[1]{\origin{#1}_\mathrm{tw}}
\newcommand{\twdestination}[1]{\destination{#1}_\mathrm{tws}}
\newcommand{\spotmarket}[1]{s_#1}
\newcommand{\pricekm}{\kappa}
\newcommand{\servicetime}{\sigma}
\newcommand{\minkm}{\mu}
\newcommand{\consecutivedriving}{\tau_{\rm n}}
\newcommand{\breaktime}{\tau_{\rm b}}
\newcommand{\sundaybreak}{\tau_{\rm s}}
\newcommand{\averagespeed}{\nu}
\newcommand{\traveltimearc}[1]{t_{#1}}
\newcommand{\traveldistancearc}[1]{d_{#1}}
\newcommand{\traveltime}[2]{\traveltimearc{#1,#2}}
\newcommand{\traveldistance}[2]{\traveldistancearc{#1,#2}}
\newcommand\nobrkhyph{\mbox{-}}
\newcommand{\linprog}[3]{%
 \begin{subequations}\allowdisplaybreaks
  \begin{alignat}{2}%
   \text{#1} \quad & #2 \\%
   \text{s.\,t.} \quad #3%
  \end{alignat}
 \end{subequations}
}
\newcommand{\constraint}[4]{ & #1 \, #2 \, #3 & \quad & \text{#4} }
\let\var\boldsymbol
\newcommand{\bin}{\mathbb B}
\newcommand{\nat}{\mathbb N}
\newcommand{\target}{\mathrm{target}}
\newcommand{\source}{\mathrm{source}}
\newcommand{\twoconds}[2]{\substack{#1\\#2}}
\newcommand{\removalLimitPercentage}{\xi}
\newcommand{\removalLimitAbsolute}{\psi}
\newcommand{\maxiter}{m}
\newcommand{\updateOperatorProbability}{n}

\journal{Journal of \LaTeX\ Templates}

\begin{document}

\begin{acronym}
\acro{ALNS}{adaptive LNS}
\acro{FCT}{full charter trucks}
\acro{FTL}{full truck load}
\acro{LNS}{large neighborhood search}
\acro{LTL}{less-then-truck load}
\acro{OVRP}{open vehicle routing problem}
\acro{PDPTW}{pickup and delivery problem with time windows}
\acro{SM}{spot-market}
\end{acronym}

\begin{frontmatter}

\title{Spot Market versus Full Charter Fleet:\\Decisions Support for Full Truck Load Tenders}

\author[1]{Veaceslav Ghilas}
\ead{veaceslav.ghilas@perfectart.com}
\author[2]{Ivo Hedtke\corref{cor1}}
\ead{ivo.hedtke@dbschenker.com}
\author[2]{Joachim Weise}
\ead{joachim.weise@dbschenker.com}
\author[3]{Tom Van Woensel}
\ead{t.v.woensel@tue.nl}
\cortext[cor1]{Corresponding author}
\affiliation[1]{organization={Perfect Art Ltd.}, city={London}, country={England}}
\affiliation[2]{organization={Global Data Strategy \& Analytics, Schenker AG}, addressline={Kruppstr.~4}, city={Essen}, postcode={45128}, country={Germany}}
\affiliation[3]{organization={School of Industrial Engineering, Eindhoven University of Technology}, country={the Netherlands}}

\begin{abstract}
This paper presents an approach to help business decision-makers gain market share by providing competitive tender offers for \ac{FTL} services.
In particular, we compare operating a fleet of \ac{FCT}, using \ac{SM} capacity and a mixture of both options against each other.
A Pickup and Delivery Problem is modeled, and solved using an Adaptive Large Neighborhood Search heuristic.
Computational results indicate strong service benefits combining \ac{FCT} and \ac{SM} usage.
Numerical experiments are presented in detail to support the findings.
Additionally, a real-life case study originating from DB Schenker is presented.
\acresetall
\end{abstract}


\end{frontmatter}

\section{Introduction}

\noindent
\Ac{FTL} transportation gained limited attention from the scientific community, compared to \ac{LTL} or express deliveries \cite{wieberneit}.
However, there is a multitude of aspects related to \ac{FTL} transportation, which have a significant financial impact and are worth researching, \eg, fleet sizing, dynamic routing, or \ac{FTL} solution design.

Freight forwarders regularly compete at tenders organized by large shippers who need tailored full-load logistics solutions.
To gain more market share, the forwarders need to provide cost-efficient full-load solutions, which are relatively complex to analyze.
The core trade-off when designing such solutions is the following: 
on the one hand, the solution needs to be economical and avoid overexposure to prohibitively high fixed costs; 
on the other hand, the shipper needs to be provided with sufficient capacity at nearly any moment, despite volume development. 

Consequently, a freight forwarder has two extreme options to organize capacity. 
One would be to cover all shipping demand exclusively through flexible capacity sourced from \ac{SM}. 
Usually, \ac{SM}~can provide service at short notice. 
This flexibility, however, comes at a relatively high market price and can not be taken for granted in special market situations, \eg, during peak seasons or around public holidays. 

Alternatively, the forwarder could invest into a fleet of \ac{FCT} dimensioned according to the (assumed) peak demand. The main idea is the following: rather than providing a customer an offer based on \ac{SM}~prices per trip or trade lane, the forwarder would instead calculate the costs of chartering and managing a fleet of trucks on a long-term basis, \eg, several months or years. When sizing this fleet properly and dispatching it efficiently, the offer involving \ac{FCT} may eventually be more competitive than an \ac{SM}~solution.

Finally, an effective mix of both solutions mentioned above may, at least in theory, provide an even more cost-efficient solution than both individual options introduced above.

The problem has similarities with a multitude of problems extensively studied by researchers over the years. On the one hand, it boils down to the \ac{PDPTW} \cite{sav_sol1995}, in particular to the \ac{FTL} variant of the \ac{PDPTW}. 

A recent study on \ac{FTL}-\ac{PDPTW} can be found in \cite{SOARES2019174}. The authors propose a meta-heuristic to tackle the problem with synchronization constraints. The approach is validated in a case study on the biomass logistics industry. In a related study, \cite{XUE2021500} develop a column generation heuristic algorithm to solve large problem instances for \ac{FTL} routing problem with multiple shifts. The described approach outperforms the previous attempts in terms of computational time and solution quality. Due to the lengthy planning horizons considered, driving regulation constraints become crucial in achieving realistic solutions. \cite{asvinGoel2009, Goel:2010:TDS:1886581.1886582, GOEL2018144} study the driving regulations and modeling approaches in the vehicle routing context. \cite{GOEL2013a} describe a heuristic algorithm designed to solve the vehicle scheduling problem and analyze different driving regulations around the world. In a related study, \cite{GOEL2016a} propose an exact approach to tackle the problem to optimality. In addition, considering the assumption that the trucks do not have pre-defined start and end depots leads us to a particular case of the well-known \ac{OVRP} \cite{palgrave.jors.2602143}, where vehicles do not have to return to a depot. A recent study on \ac{OVRP} can be found in \cite{doi:10.1080/00207543.2019.1572929}. The authors propose a hybrid adaptive large neighborhood search approach to solve large-scale problems for \ac{OVRP} and manage to find new best-known solutions for specific instances. Furthermore, various routing problems with practical side constraints, such as \ac{FTL} property and driving regulations, have been studied by scientists in the last decades. The interested reader is referred to \cite{RePEc:eee:ejores:v:231:y:2013:i:1:p:1-21, Braekers2016TheVR, lit_review2016} for recent surveys on different types of routing problems. 

The orienteering problem \cite{doi:10.1057/palgrave.jors.2602603} is another related well-known problem. In this context, each resource (\ie, truck) aims at finding a minimum cost (or the most profitable) path, starting at a pickup location and ending at a delivery location. Expensive requests, in terms of \ac{FCT}~cost, would then be serviced by \ac{SM}. In this sense, the aim is to find the right trade-off between~\ac{FCT} and \ac{SM}~costs. \cite{orient_pdp} propose a general variable neighborhood search algorithm to solve large-size orienteering problems and show the efficiency of the developed approach.  In addition, \cite{sun} investigate a time-dependent orienteering problem and propose a tailored adaptive large neighborhood search. The authors show the benefits of considering time-dependency during the planning process to obtain efficient and reliable solutions.

The contributions of the paper at hand are as follows:
\begin{enumerate}
\item We analyse the needed decision support for \ac{FTL} transportation and the related tendering process.
\item We model the described problem as a Pickup and Delivery Problem and propose a tailored adaptive large neighborhood search (ALNS) to solve the underlying routing problem.
\item Numerical results show that the mix of \ac{SM} and \ac{FCT} may lead to cheaper FTL solutions compared to traditional pure \ac{SM} strategy, having the transformed VRP Gehring \& Homberger instances as benchmarks.
\item We also present a practical case application to the tender support for \ac{FTL} business based on a real-life setting originating from DB Schenker. DB Schenker is a freight forwarder that supports industry and trade in the global exchange of goods by land transport, worldwide air and ocean freight, contract logistics and supply chain management.
\end{enumerate}

This paper is organized as follows. In Section~\ref{section:problemdescription} we formally define the considered mathematical problem, we present the related literature, and discuss the transformation of the business scenario into the formal problem. In Section~\ref{section:solutionapproach} we present a single-solution meta-heuristic tailored to solve the problem at hand. Finally, we discuss the experiments on transformed literature instances and real-world data sets involving two customers and the obtained solutions in Section~\ref{section:computationalexperimentsandresults}.

\section{Problem Description and Mathematical Model}
\label{section:problemdescription}

\noindent
We consider a set~$\requests$ of transportation requests between a finite set of locations~$\locations$.
Each request~$r\in \requests$ has an origin~$\origin{r}\in \locations$, a destination~$\destination{r} \in \locations$, 
as well as a time window~$\tworigin{r}$ at the origin and a set of time windows~$\twdestination{r}$ at the destination.
A time window is defined by two absolute points in time.
Every request must be picked up at its origin within its origin time window, transported directly to its destination and be delivered within one of its destination time windows.
The travel-time and travel-distance between two locations~$\ell,\ell'$ are given by~$\traveltime{\ell}{\ell'}$ and~$\traveldistance{\ell}{\ell'}$.

Each transportation request~$r$ can be outsourced at a cost of $\spotmarket{r}$.
The remaining requests must be served by a set of vehicles at a cost of
$\pricekm$ per unit of driven distance.

The task at hand is to find a cost-optimal assignment of all requests to the options of outsourcing it or serving it by a vehicle. Note that the set of vehicles is determined as part of the task.

Serving the non-outsourced requests by vehicles is subject to the following constraints:
\begin{itemize}
\item 
Each vehicle can serve at most one request at any given point in time.
\item 
Each vehicle starts at the origin of its first request, ends at the destination of its last request,
and has to drive a distance of at least $\minkm$.
\item
Each (un\nobrkhyph)loading operation requires a time of~$\servicetime$.
The complete (un\nobrkhyph)loading operation has to be executed within one of the respective time windows.
Hence, if a vehicle arrives outside the time windows, it has to wait for the next one.
\item 
Each vehicle has to fulfill the following driving time regulations:
A shift-break is a break of duration of at least~$\breaktime$ and a Sunday-break is a break of duration of at least~$\sundaybreak \geq \breaktime$.
After a duration of~$\consecutivedriving$ cumulative driving without shift-breaks, a shift-break is required.
\mbox{(Un\nobrkhyph)}Loading does not count as break time, whereas waiting time does.
Every Sunday, a Sunday-break is required.
\end{itemize}

A mathematical model combining classical arc-based vehicle routing, the minimum driving distance constraint from \cite[Section~III]{Kara2011}, and the labeling technique from \cite{GoelGruhn2006} for the driving time regulations is constructed in three phases:
First, we define a graph as basis for an arc-based vehicle routing model.
Then, a model respecting everything except for the timing constraints is built.
Finally, the timing constraints (time windows, simplified driving regulation) are added.

Variables are denoted by bold lower-case letters.
Instance independent parameters are denoted by greek letters.
We write $\bin :=\{0,1\}$ and $\nat :=\{0,1,2,\ldots\}$.
By $M$ we denote a sufficiently large constant.
An overview of the variables and parameters can be found in Table~\ref{tab:ModelElements}.

\begin{table}
\centering\footnotesize
\begin{tabular}{ll}
\hline
Variable & Meaning\\\hline
$\var x_a$        & vehicle is traveling along arc~$a$ (=1) or not (=0)\\
$\var y_{n\to m}$ & distance from origin to~$m$ when a vehicle travels along $n\to m$\\
$\var l_n=(l_{n,1},l_{n,2})$ & state of the driver at node~$n$ (arrival time, nonstop driving time)\\\hline\hline
Parameter & Meaning\\\hline
$\traveldistancearc a$ & distance of arc~$a$\\
$\pricekm$ & cost of driving per unit distance\\
$\minkm$ & minimum driving distance of a vehicle\\
$M$ & sufficiently large constant\\
$\spotmarket{r}$ & cost of out-sourcing request~$r$ to the spot-market\\\hline
\end{tabular}
\caption{Parameters and variables in the model.}
\label{tab:ModelElements}	
\end{table}

Let $G:=(N,A)$ denote the directed graph constructed as follows: The set~$N$ of nodes consists of
\[
N:=\{(\origin r, r) : r\in R\} \cup \{(\destination r, r) : r\in R\} \cup \{n_0,n_\infty\}
\]
all origin- and destination locations together with two artificial nodes~$n_0$ and~$n_\infty$ denoting the start and end of all tours.
We write~$N_R:=N\setminus\{n_0, n_\infty\}$.
For a request~$r$ we define~$r^s$ as the start of the time window~$\tworigin r$ and~$r^e$ as its end.
The set~$A$ of arcs represents the vehicles realizing the non-outsourced request, traveling empty between locations or waiting at a location.
It consists of the arcs serving requests $r\in R$, the tour start, the tour end, and waiting/empty travel:
\begin{align*}
A:={}&\{(\origin r,r) \to (\destination r,r) : r\in R\}\\
&\cup \{n_0 \to (\origin r,r) : r\in R\}\\
&\cup \{(\destination r,r) \to n_\infty : r \in R\}\\
&\cup \left\{(\destination{r_1},r_1) \to (\origin{r_2},r_2) : r_1\neq r_2 \in R, r_s + \traveltime{\origin{r_1}}{\destination{r_1}} +\traveltime{\destination{r_1}}{\origin{r_2}} + 3\servicetime \leq r^e_2\right\}.
\end{align*}
For an arc~$a=n\to m$ we define $\target (a):=m$ and $\source(a):=n$ as the target and source of~$a$.
We extend the distances $\traveldistance{\ell}{k}$ and travel times $\traveltime{\ell}{k}$ between locations to $\traveldistancearc a$ and $\traveltimearc a$ on arcs $a\in A$ by
\[
\traveldistancearc a := \begin{cases}	
0 & \source (a) = n_0,\\
0 & \target(a) = n_\infty,\\
\traveldistance{\ell}{k} & a=(\ell, r) \to (k, r')
\end{cases}
\quad\text{and}\quad
\traveltimearc a := \begin{cases}	
0 & \source (a) = n_0,\\
0 & \target(a) = n_\infty,\\
\traveltime{\ell}{k} & a=(\ell, r) \to (k, r').
\end{cases}
\]

Let $\var x_a \in \bin$ denote whether a vehicle is traveling along arc~$a\in A$ ($\var x_a = 1$) or not ($\var x_a = 0$).
Note that for an arc $a=(\origin r, r)\to(\destination r, r)$ this denotes whether the request~$r$ is served by vehicles ($\var x_a = 1$) or out-sourced ($\var x_a = 0$).
As introduced in \cite[Section~III]{Kara2011} we define $\var y_{n\to m}\in \nat$ as the total distance from the origin $n_0$ to node $m\in N$ traveled by a vehicle when it goes along $n\to m \in A$.
We focus purely on the routing and the minimum distance per vehicle constraint:
\linprog{min}{%
 \label{m:obj}%
 \sum_{a\in A}  \pricekm \traveldistancearc a \var x_a +  \sum_{r\in R} s_r (1 - \var x_{(\origin r,r) \to (\destination r,r)})}%
{%
\constraint{%
 \label{m:kirchhoff1}%
 \sum_{\twoconds{a \in A}{\target(a)=(\origin r, r)}}\var x_a}{=}{\var x_{(\origin r, r) \to (\destination r, r)}}{$\forall r \in R$}\\
\constraint{%
 \label{m:kirchhoff2}%
 \sum_{\twoconds{a\in A}{\source(a) = (\destination r,r)}}\var x_a}{=}{\var x_{(\origin r,r) \to (\destination r,r)}}{$\forall r \in R$}\\
\constraint{%
 \label{m:drivupdate}%
 \sum_{\twoconds{a\in A}{\source(a) = n}}\var y_a}{=}{\sum_{\mathclap{\twoconds{a\in A}{\source(a) = n}}}d_a \var x_a + \sum_{\mathclap{\twoconds{a\in A}{\target(a) = n}}}\var y_a}{$\forall n \in N_R$}\\
\constraint{%
 \label{m:n0}%
 \var y_{n_0 \to (\origin r, r)}}{=}{0}{$\forall r\in R$}\\
\constraint{%
 \label{m:bigm}%
 \var y_a}{\leq}{M \var x_a}{$\forall a\in A$}\\
\constraint{%
 \label{m:mindriving}%
 \var y_a}{\geq}{\minkm \var x_a}{$\forall a\in A$: $\target(a)=n_\infty$}\\
\constraint{%
 \label{m:vars}%
 \var x_a \in \mathbb B}{,}{\var y_a\in \mathbb R^+}{$\forall a \in A$}
}

Finally, we consider the simplified driving time regulations. Inspired by \cite{GoelGruhn2006}, we define for all nodes~$n\in N_R$ a label: 
\[
\var l_n=\begin{pmatrix}
l_{n,1}\\
l_{n,2}
\end{pmatrix}
=
\begin{pmatrix}
\text{arrival time}\\
\text{nonstop driving time}
\end{pmatrix}
\]
to represent the state of the driver at the node.
The vehicle can start the service at node~$n\in N_R$ at time $l_{n,1}$ and can depart from~$n$ at time~$l_{n,1}+\servicetime$.
It may drive~$\consecutivedriving - l_{n,2}$ before the next break.

As shown in \cite[Section~V]{GoelGruhn2006}, it is possible to pre-compute the set~$\mathcal L_m(\var l_n)$ of potential labels for a vehicle that is supposed to travel from node~$n\in N_R$ with label~$\var l_n$ to a node~$m\in N_R$.
This~takes
\begin{itemize}
\item the simplified driving regulation,
\item the time window(s) of the requests~$r\in R$ at origin and destination, and
\item the Sunday-break~$\sundaybreak$
\end{itemize}
into account.
Let $l^r:=(r^s, 0)^\top$ denote the label of a vehicle when the request~$r$ is the first request within the tour.

The extended model taking time windows and driving time regulation into account is then
\linprog{min}{\eqref{m:obj}\notag}{%
\constraint{\eqref{m:kirchhoff1}, \eqref{m:kirchhoff2}, \eqref{m:drivupdate}, \eqref{m:n0}, \eqref{m:bigm}, \eqref{m:mindriving}, \eqref{m:vars}}{}{}{\notag}\\
\constraint{\var x_{n_0\to (\origin r,r)} = 1}{\Longrightarrow}{\var l_{(\origin r,r)} = l^r}{$\forall r \in R$}\\
\constraint{\var x_{n\to m} = 1}{\Longrightarrow}{\var l_m\in \mathcal L_m(\var l_n)}{$\forall n\to m \in A$: $n,m\in N_R$}
}

\section{Solution Approach}
\label{section:solutionapproach}

\noindent
The algorithm presented in this section is a single-solution meta-heuristic based on \ac{LNS}. 
The general idea originates in \cite{Shaw1998}, and has been gaining significant popularity in the recent past. 
In particular, it has proven to be highly effective in tackling vehicle routing problems and provides a good trade-off between solution quality and computational time.
Note that we will use \emph{shipment} as a synonym for request in the following section.

The main workflow of the algorithm is as follows: given an initial solution, iteratively modify it until a stopping criterion is reached. 
Two basic operators in the \ac{LNS} framework are \emph{removal} and \emph{insertion} procedures. 
In the current context, \ie, the pickup and delivery problem, the removal method un-plans a percentage~$\removalLimitPercentage$ of already planned requests (pickup and delivery vertices associated with a request) from the solution, and the insertion method would re-plan these requests, by using, \eg, an insertion heuristic.
Note that for large instances, removing, \eg, 20\,\% of the shipments might still lead to a large optimization sub-problem, hence we limit it by an absolute upper bound~$\removalLimitAbsolute$ on the number of requests to remove.
At every iteration, the modified solution is compared to the best solution found, and promoted to the new best solution if proved to be better.
It is worth noting that a solution might contain unplanned requests, which are assumed to be serviced by SM, and thus incur additional costs.
The original framework is extended to consider an adaptive mechanism \cite{Pisinger2007ALNS}, \ie, \ac{ALNS}, where user-defined removal and insertion methods are selected based on their performance during the search.

The high-level pseudo-code of \ac{ALNS} is presented in Algorithm~\ref{algo:alns}.
The interested reader is referred to \cite{Pisinger2007ALNS, ropke2006, Ghilas:2016:ALN:2927990.2928062, GRIMAULT20171} for more details about applying \ac{ALNS} to a broad range of VRPs.

\begin{algorithm}[t]\footnotesize
	\caption{ALNS}\label{algo:alns}
	\begin{algorithmic}
		\State $s$ $\gets$ \texttt{generateInitialSolution()}
		\State $s_\mathrm{best}$ $\gets$ $s$
		\While{!\textit{stoppingCriterion}}
			\State $r^*$ $\gets$ \texttt{chooseRemovalOperator()}
			\State $i^*$ $\gets$ \texttt{chooseInsertionOperator()}
			\State $s'$ $\gets$ \texttt{partiallyDestroySolution($r^*$, $s$)}
			\State $s_\mathrm{new}$ $\gets$ \texttt{repairSolution($i^*$, $s'$)}
			\If {\texttt{accept($s_\mathrm{new}$)}}
				\State $s \gets s_\mathrm{new}$
			\EndIf
			\If {$s_\mathrm{new}$ \textless{} $s_\mathrm{best}$}
			    \State $s_\mathrm{best}$ $\gets$ $s_\mathrm{new}$
			\EndIf
			\State \texttt{updateOperatorProbabilities()}
		\EndWhile
		\Return $s_\mathrm{best}$
	\end{algorithmic}
\end{algorithm}

The starting solution~$s$ is constructed using a greedy algorithm.
Its basic idea is to insert an unplanned request in the feasible position which increases the objective function value the least.
Note that the initial solution out-sources all shipments, nothing is planned on vehicles.
A request is given to a vehicle only if the insertion cost is smaller than the corresponding out-sourcing cost.
While the stopping criterion is not reached, \eg, the maximum number of iterations~$\maxiter$, the algorithm randomly selects one removal and one insertion operator, and applies them to the current solution $s$, thus generating a new solution $s_\mathrm{new}$.
The method \texttt{accept} verifies whether the newly-built solution is accepted.
Simulated annealing acceptance criteria are used in the current implementation, similarly to \cite{Pisinger2007ALNS, ropke2006}, so that worse solutions may also get accepted.
If the newly-built solution is accepted, the current solution~$s$ is overwritten.
Finally, the best solution found is returned.

Note that the probability of an operator being chosen is dynamically updated every~$\updateOperatorProbability$th iteration.
The better the performance of the operator, the higher its chance of being chosen.
Eventually, the algorithm will converge to using only good-performing operators, see \cite{Pisinger2007ALNS, ropke2006}.

It is important to note that the number of vehicles is unlimited, however at the end of each \ac{ALNS} iteration, certain number of requests may remain unplanned.
This is mainly due to the fact that out-sourcing may be cheaper if a significant empty travel would be induced by servicing a specific request using a vehicle, or if simply the \ac{SM} rate is cheaper than the \ac{FCT} rate.

\subsection{Removal Operators}

\noindent
Several removal operators are used in the current implementation, as described below. 

\begin{description}
\item[RRR:]
Random Route Removal randomly selects a route and removes it from the solution.
All routes have the same probability of being selected;
\item[TRR:]
Time-based Route Removal is similar to RRR, however, the probability of a route being chosen depends on the total travel time of the corresponding route in the current solution.
In other words, longer total travel time leads to a higher probability of being selected;
\item[SRR:]
Stop-based Route Removal is similar to TRR, however, smaller number of shipments planned within a route leads to a higher probability of being selected.
The intuition behind this is that the fewer requests are planned in a route, the easier it is to re-plan them into other routes, thus avoiding this route in the solution.
\item[RSR:]
Random Shipment Removal randomly selects a set of requests to be removed from the current solution.
This operator helps in terms of search diversification. 
\item[TSR:]
Time-based Shipment Removal is similar to RSR, however, the probabilities of being selected depend on the incurred driving times.
In particular, the probability of selecting a shipment is higher if, in a given solution, total driving time to its pickup location, and from its delivery location, is longer.
\item[SR:]
Shaw Removal removes a set of shipments similar to each other (\cite{Shaw1998}).
The similarity is defined by the distance between pickup locations and delivery locations of each pair of shipments.
In addition, the time windows of the shipments are part of the similarity function.
Two variants of SR are used: (\emph{i}) with distance and time windows as similarity criteria, and (\emph{ii}) only time windows as similarity criterion.
\end{description}

\subsection{Insertion Operators}

\noindent
We provide more details about the insertion operators used as follows.

\begin{description}
\item[Classical greedy:]
identifies in each iteration the request (of the set of unplanned requests and the set of (partial) routes) which incurs the lowest insertion cost and inserts it in its best feasible position.
It repeats this operation until no unplanned request exists or no insertion cost (\ie, \ac{FCT} cost) is lower then the corresponding \ac{SM} cost.
\item[Regret insertion:]
finds the shipment which incurs the maximum \emph{regret} if not inserted in its cheapest feasible position at every iteration.
Let $c_1$ represent the cost of inserting the shipments in the route with the cheapest feasible position within the \ac{FCT} routes, $c_2$ -- in another route with the second cheapest position, $c_3$ -- third, etc.
Then, the \emph{regret function} can be defined as $c_2-c_1$, known as $2$-regret. For more details on regret insertion, please refer to \cite{potvin1993}.
In order to take into account more information when deciding which request to insert next, the regret function is generalized, known as $k$-regret, considering multiple routes.
For example, it is possible to look at the three cheapest insertion positions. \ie, $\sum_{i=2}^k(c_k - c_1)$, where $k = 3$ is the number of look-ahead insertion positions to take into account (\ie, $3$-regret). Note that several regret operators (with different $k$ values) can be used within the \ac{ALNS}.
\end{description}

\begin{algorithm}[t]\footnotesize
	\caption{Insertion procedure}\label{algo:insert}
	\begin{algorithmic}
		\Require $S_\mathrm{in}=\mathit{unplannedShipments}$, $P=\mathit{partialRoutes}$, $\spotmarket r$ (out-sourcing costs)
		\State $S_\mathrm{out} \gets \emptyset$
		\While{$S_\mathrm{in} \neq \emptyset$}
			\For{$r \in S_\mathrm{in}$, $p\in P$}
				\State $c_{rp}$ $\gets$ \texttt{computeCostOfInsertingRequestInRoute($r$, $p$)}
			\EndFor
			\State $(r^*, p^*)$ $\gets$ \texttt{shipmentToInsertNext($c$)}
			\If{$c_{r^*p^*} < \spotmarket{{r^*}}$}
			    \State \texttt{insertRequestInRoute($r^*$, $p^*$)}
			\Else
			    \State $p_\mathrm{new} \gets \emptyset$
			    \State \texttt{insertRequestInRoute($r^*$, $p_\mathrm{new}$)}
			    \If {\texttt{allRoutesAreWellUtilized($P$)} \textbf{and} \texttt{cost($p_\mathrm{new}$)}${} \leq \spotmarket{{r^*}}$}
			        \State $P=P \cup \{p_\mathrm{new}\}$
			    \Else
			        \State $S_\mathrm{out} \gets S_\mathrm{out} \cup \{r^*\}$
			    \EndIf
			\EndIf
			\State $S_\mathrm{in} \gets S_\mathrm{in} \setminus r^*$
		\EndWhile
		\Return $S_\mathrm{out}$
	\end{algorithmic}
\end{algorithm}

The general framework of an insertion operator is shown in Algorithm~\ref{algo:insert}:
Given a set~$S_\mathrm{in}$ of unplanned shipments and a set~$P$ of partial routes, iteratively find the cost of inserting a shipment into a route that incurs the most beneficial change of the given \emph{cost function}. 
If the insertion cost is cheaper than the corresponding out-sourcing cost, the insertion is performed.

Otherwise, a new trip is created only if no feasible insertion is found in the existing vehicle trips, and the corresponding out-sourcing cost is higher or equal than using a new vehicle trip. In addition, a new trip is created only if all vehicles used are utilized w.r.t. distance traveled at least~$\minkm$ per planning horizon. 

Otherwise, the request is placed into the out-sourcing bank~$S_\mathrm{out}$, \ie, the set of the requests which are out-sourced.

The algorithm is assumed to start with one vehicle available. As soon as the vehicles used are well-utilized as aforementioned, the algorithm makes an additional vehicle available.

We compute an insertion cost matrix~$c_{rp}$ for each unplanned shipment and each existing route. 
The method \texttt{shipmentToInsertNext} returns the shipment with the least incurred insertion cost, along with the route and the corresponding position within the route.
Here, either classical greedy or a $k$-regret insertion is used.
The selected shipment is then evaluated, \ie, if the incurred insertion cost is cheaper compared to out-sourcing, it is inserted in the chosen route (in its cheapest position).
The request is then removed from~$S_\mathrm{in}$. 
Note that in the event that no feasible insertion is found, the shipment $r$ with the least cost per distance unit is selected (\ie, $\spotmarket{r}/d_{s_{o},s_{d}}$).
Finally, the insertion procedure returns the out-sourcing bank~$S_\mathrm{out}$.

\subsection{Constraints}

\noindent
The feasibility of (intermediate) solutions is enforced at any time during the run of the algorithm. 
To recap, the following constraints need to be satisfied during cost matrix computation: (\emph{i}) every request is out-sourced or must be served by a vehicle, (\emph{ii}) time window(s) at origin and destination of the requests must be respected, and (\emph{iii}) breaks between working shifts as well as (\emph{iv}) Sunday-breaks must be taken into account.
Constraints (\emph{ii}), (\emph{iii}), and (\emph{iv}) are enforced only if the requests are served by a vehicle.
As aforementioned, out-sourced requests are assumed to satisfy all the considered constraints.

Constraints (\emph{i}) and (\emph{ii}) are straightforward to implement.
Many researchers have investigated ways to consider these constraints in an efficient manner, \ie, by using auxiliary data structures, \eg, \cite{savelsbergh1990,campbell2004}.
However, when combined with constraints (\emph{iii}) and (\emph{iv}), it is not trivial to efficiently implement them within \ac{ALNS}.

\section{Computational Experiments}
\label{section:computationalexperimentsandresults}
\noindent
First, we describe the assumptions made during the transformation of Gehring \& Homberger VRP instances, along with their corresponding results.
The second part presents a case study at DB Schenker and describes the instance characteristics of the input data sets and costs.
In both sections, the costs for out-sourcing no request and out-sourcing all requests are computed for comparison.
Finally, the solutions that assign all requests to outsourcing it or serving it by a vehicle are discussed.

The ALNS is implemented in C++11 and all experiments are run on an Intel~Core~i7-8750H machine @2.2GHz, with 16~GB DDR4-RAM @2.4GHz.

\paragraph{Algorithm Parameters}
The algorithm contains various parameters that need to be set. 
Some of them are business-related, others are heuristics technical in nature. The business-related parameters, along with corresponding explanation and values, are shown in Table~\ref{tab:business_params}.

\begin{table}
\centering\footnotesize
\begin{tabular}{llrl}
\hline
Parameter & Description & Value & Unit \\
\hline
$\consecutivedriving$ & duration of cumulative driving without shift-breaks & 450 & min \\
$\breaktime$ & minimum duration of a break & 990 & min \\
$\sundaybreak$ & minimum duration of a Sunday-break & 1,320 & min \\
$\servicetime$ & duration of an (un-)loading operation & 120 & min \\
$\averagespeed$ & average speed of vehicles & 70 & km/h \\\hline
\end{tabular}
\caption{Business-specific parameters.\label{tab:business_params}}
\end{table}

After performing extensive computational experiments, we found that the best performance can be achieved by applying the removal operators RRR, SRR, SR, TSR, and RSR and the insertion operators greedy, $4$-regret, $5$-regret, and $6$-regret, as well as the technical parameter settings of Table~\ref{tab:tech_params}.

\begin{table}
\centering\footnotesize
\begin{tabular}{llr}
\hline
Parameter & Description & Value\\
\hline
$\maxiter$ & maximum number of \ac{ALNS} iterations & 25,000 \\
$\removalLimitAbsolute$ & maximum absolute number of shipments to remove & 100 \\
$\removalLimitPercentage$ & maximum relative number of shipments to remove  & 35\,\% \\
$\updateOperatorProbability$ & every $\updateOperatorProbability$ iterations update operators weights & 200 \\
\hline
\end{tabular}
\caption{Technical parameters.\label{tab:tech_params}}
\end{table}

Furthermore, we apply simulated annealing (SA) acceptance criteria in our ALNS framework, as inspired by \cite{ropke2006}. We considered a similar parameter setup as in \cite{Ghilas:2016:ALN:2927990.2928062}, as it proved to be beneficial for solution quality.

Note that this presented parameter setup is used throughout this section.

\subsection{Transformed Gehring \& Homberger Instances}

\noindent
In this section, we present the computational results obtained by solving the transformed well-known Gehring \& Homberger VRP instances \cite{GehringHomberger1999}. To convert the literature instances, we made the following assumptions:
\begin{itemize}
    \item Node 0 is ignored, as it corresponds to the depot;
    \item Demand, capacity and service time data is ignored;
    \item For an instance with $N$ nodes, node $n$ and $N$/2 + $n$ correspond to pickup and delivery nodes of a request;
    \item To assure that the planning horizon consists of multiple days, we multiplied the Euclidean distances and pickup start time windows by a factor $f = 6$. The resulting distances are then rounded to the closest integer;
    \item For simplicity, all locations (\ie, pickup and delivery locations) are assumed to be open between 06:00 and 18:00; \item The transformed pickup time is used to determine the ready day within the planning horizon. Then, the corresponding time window from 06:00 to 18:00 of that day is used as pickup time window;
    \item \ac{SM} rates per km are assumed as shown in Table \ref{tab:smRates_lit};
    \item \ac{FCT} rate is assumed 1.06\,EUR/km;
    \item $\mu$ is computed for each instance separately: 250\,km per day with pickups in the planning horizon.
\end{itemize}

Table \ref{results_lit} presents the results obtained after solving the transformed Gehring \& Homberger instances. In particular, the numbers show the averages over all corresponding instances for each instance class (\eg, C100). Three scenarios are computed: all requests served by an own fleet, all requests out-sourced, and a mix between those as mentioned above.

\begin{table}
\centering\footnotesize
\begin{tabular}{l|ccc}%
\hline%
distance~$\delta$ & $\delta < {}$150\,km & 150\,km${}\leq \delta <$350\,km & 350\,km${}\leq \delta$\\
\hline
cost in EUR/km & 1.75 & 1.40 & 1.15\\
\hline
\end{tabular}
\caption{Assumptions on spot-market rates.\label{tab:smRates_lit}}
\end{table}

\begin{table}
\centering\footnotesize
\begin{tabular}{@{}l|r@{~}r|r@{~}r|r@{~}r@{~}r@{~}r|r@{~}r|r@{}}
\hline
 & \multicolumn{2}{c|}{out-sourcing} & \multicolumn{2}{c|}{out-sourcing} & \multicolumn{7}{c}{mixed scenario}\\
Class & \multicolumn{2}{c|}{nothing} & \multicolumn{2}{c|}{everything} & \multicolumn{4}{c|}{own vehicles} & \multicolumn{2}{c|}{out-sourced} & $\Sigma$\\
\&  & \multicolumn{2}{c|}{} & \multicolumn{2}{c|}{} & \% & km & km & &&& \\
Nodes & km & cost & km & cost & req. & loaded & empty & cost & km & cost & cost\\
\hline
C\phantom{R} 100 &    58.0 &      64.5 &   45.0 &     54.6 &     2.2 &    1.14 &   0.04 &           1.3 &       44.0 &         53.2 &      \bf  54.5 \\
C\phantom{R} 200 &   145.0 &     154.5 &  115.0 &    136.4 &     5.2 &    6.00 &   0.18 &           6.5 &      109.0 &        129.0 &      \bf 135.6 \\
C\phantom{R} 300 &   324.0 &     343.3 &  267.0 &    310.1 &     5.2 &   17.29 &   0.41 &          18.8 &      250.0 &        290.0 &      \bf 308.8 \\
C\phantom{R} 400 &   527.0 &     558.3 &  449.0 &    518.6 &     5.6 &   29.70 &   0.74 &          32.3 &      419.0 &        484.2 &      \bf 516.5 \\
C\phantom{R} 500 &   869.0 &     921.2 &  764.0 &    881.0 &     7.2 &   67.31 &   1.15 &          72.6 &      697.0 &        803.4 &      \bf 876.0 \\
\hline
R\phantom{C} 100 &    58.0 &      63.3 &   44.0 &     53.0 &     2.0 &    0.90 &   0.02 &           1.0 &       43.0 &         51.9 &      \bf  52.9 \\
R\phantom{C} 200 &   154.0 &     163.6 &  120.0 &    142.2 &     2.0 &    2.35 &   0.08 &           2.6 &      118.0 &        139.4 &      \bf 141.9 \\
R\phantom{C} 300 &   342.0 &     362.7 &  272.0 &    314.9 &     1.7 &    6.06 &   0.11 &           6.5 &      266.0 &        308.0 &      \bf 314.5 \\
R\phantom{C} 400 &   607.0 &     643.0 &  497.0 &    573.8 &     2.3 &   15.26 &   0.32 &          16.5 &      482.0 &        556.2 &      \bf 572.7 \\
R\phantom{C} 500 &   960.0 &    1017.9 &  791.0 &    910.9 &     3.2 &   34.09 &   0.49 &          36.7 &      757.0 &        871.7 &      \bf 908.4 \\
\hline
          RC 100 &    59.0 &      63.1 &   44.0 &     53.4 &     2.1 &    1.20 &   0.03 &           1.3 &       43.0 &         52.0 &      \bf  53.3 \\
          RC 200 &   155.0 &     164.5 &  126.0 &    147.7 &     2.4 &    3.55 &   0.08 &           3.8 &      122.0 &        143.5 &      \bf 147.4 \\
          RC 300 &   338.0 &     357.8 &  262.0 &    305.2 &     2.0 &    6.35 &   0.19 &           6.9 &      256.0 &        297.8 &      \bf 304.7 \\
          RC 400 &   622.0 &     659.2 &  500.0 &    577.8 &     2.5 &   17.82 &   0.43 &          19.3 &      482.0 &        557.3 &      \bf 576.6 \\
          RC 500 &   980.0 &    1039.2 &  806.0 &    929.5 &     2.8 &   30.08 &   0.69 &          32.6 &      776.0 &        894.8 &      \bf 927.5 \\
\hline
\end{tabular}
\caption{Costs (average: grouped by class and number of nodes) for the mixed scenarios. Distances are given in $1{,}000$\,km. \%req.\ denotes the percentage of requests. Costs are given in $1{,}000$\,EUR. The highlighted cost is better.}
\label{results_lit}
\end{table}

Not surprisingly, the results indicate that considering the mix of \ac{FCT} and \ac{SM} leads to best results w.r.t. operating costs.
In particular, on average 0.5\,\% savings can be achieved for the clustered (C) instances by serving approx. 5\,\% of the requests using own vehicles compared to out-sourcing everything.
On the other hand, an average of 0.2\,\% cost savings can be achieved in the mixed scenario for random (R) and randomly-clustered (RC) instances by serving approx.\ 2\,\% of the requests using own vehicles. 

For the C instances, more savings can be achieved compared to R and RC instances.
The main reason is that requests that are compatible from a temporal point of view are more compatible from a geographical point of view.
In other words, the chances that a pickup location of a request is close to a delivery location of another request are higher in C instances. Hence less empty driving can be achieved.

Note that serving all requests by own vehicles is the most costly out of all scenarios considered. A minimum-driving-distance constraint is enforced, and the requests are not perfectly compatible in terms of timing and spatial aspects.

\subsection{DB Schenker Case Study}
\label{section:computationalexperimentsandresultsschenker}

\noindent
DB Schenker is one of the key players in the global logistics sector. 
Founded in 1872 by Gottfried Schenker in Vienna, Austria, Schenker~\&~Co.\ began its business by consolidating rail consignments from Paris, France, to Vienna, Austria. 
DB Schenker is a freight forwarder that supports industry and trade in the global exchange of goods by land transport, worldwide air and ocean freight, contract logistics, and supply chain management. 
With more than 76,000 employees working in approx.\ 2,000 locations around the world, DB Schenker is a leader in its industry.

DB~Schenker's land transport in Europe covers 36 countries and offers a variety of products and services. 
One of them is \textit{direct} product for large loads, \eg, \ac{FTL}s that are transported directly from consignor to consignee.

Large shippers frequently conduct tenders that require logistics providers to tailor dedicated full-load solutions.
DB Schenker participates in such tenders via sales department representatives, who need to prepare the business offers.
To come up with competitive solutions, different scenarios need to be analyzed, such as:
operating a fleet of \ac{FCT}, using \ac{SM} capacity, as well as a mix of both options mentioned above. 

The approach presented in this paper aims at helping sales departments get more insights about various scenarios, given \eg, historical/forecasted \ac{FTL} transportation demands. 
As a result, sales teams can become more competitive during tenders.

Two of DB Schenker's customers conducted a tender, and the sales department was given the following problem (per customer):

A set of plants is operating on a weekly schedule with working hours per day.
For an entire month, all full-load shipments that have to be transported between the plants are given.
Every shipment has a pickup date and must be transported directly (no consolidation, no pickup on a later date, etc.) to its destination.

For small instances, an offer can be created by assuming that everything is given to the spot market or that the volume can be included in DB Schenker's internal transportation network.
However, for large instances as given by the two customers (with up to 12,427 full-load shipments per month), none of the approaches has enough capacity, and none of them would enable the sales department to present a competitive offer.
An option with enough capacity that can be operated at a competitive cost is the combination of its own dedicated fleet of trucks and outsourcing to the spot market.

To determine the size of such a fleet and the cost of operating it, the above model can be applied:
The full-load shipments correspond to requests~$\requests$ that have to be transported between the locations~$\locations$ given by the plants.
The origin time window~$\tworigin r$ of a request~$r$ is given by the operating hours of the pickup plant.
The destination time windows~$\twdestination r$ is given by the operating hours of the destination plant.
Note that up to two such time windows are required, as a truck can arrive late and wait for the next day.
The two options of using the spot market and operating an own fleet correspond to outsourcing requests and serving them by a set of vehicles, respectively.

As a result should be an offer within the tender process, it is not required to create detailed schedules for the trucks of the own fleet.
A simplified driving regulations model based on shifts and including the Sunday-driving ban is sufficient.

\subsubsection{Data Description}
\label{subsection:datadescription}

\noindent
In this section, we describe the request- and cost-related data we used in the analysis.

\begin{figure}[t]
\centering
\includegraphics[trim=1in 7.2in 1in 1.1in,clip,width=\linewidth]{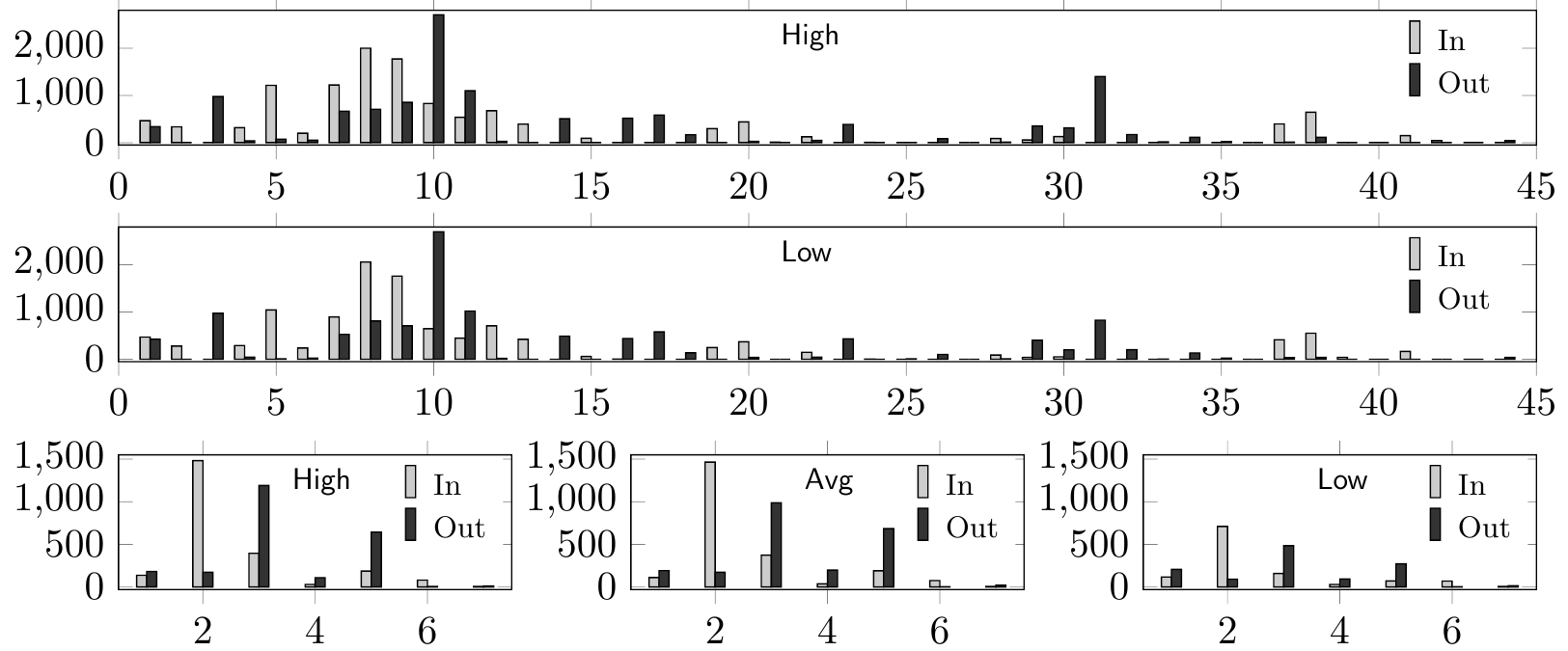}
\caption{In- and out-going requests per physical location.\label{fig:degrees} The plots show the number of in- (gray) and out-going (black) requests per location. The top and middle plot refer to Customer~1 and the remaining ones for Customer~2.}
\end{figure}

\begin{figure}[t]
\centering
\fbox{\includegraphics[width=0.48\linewidth]{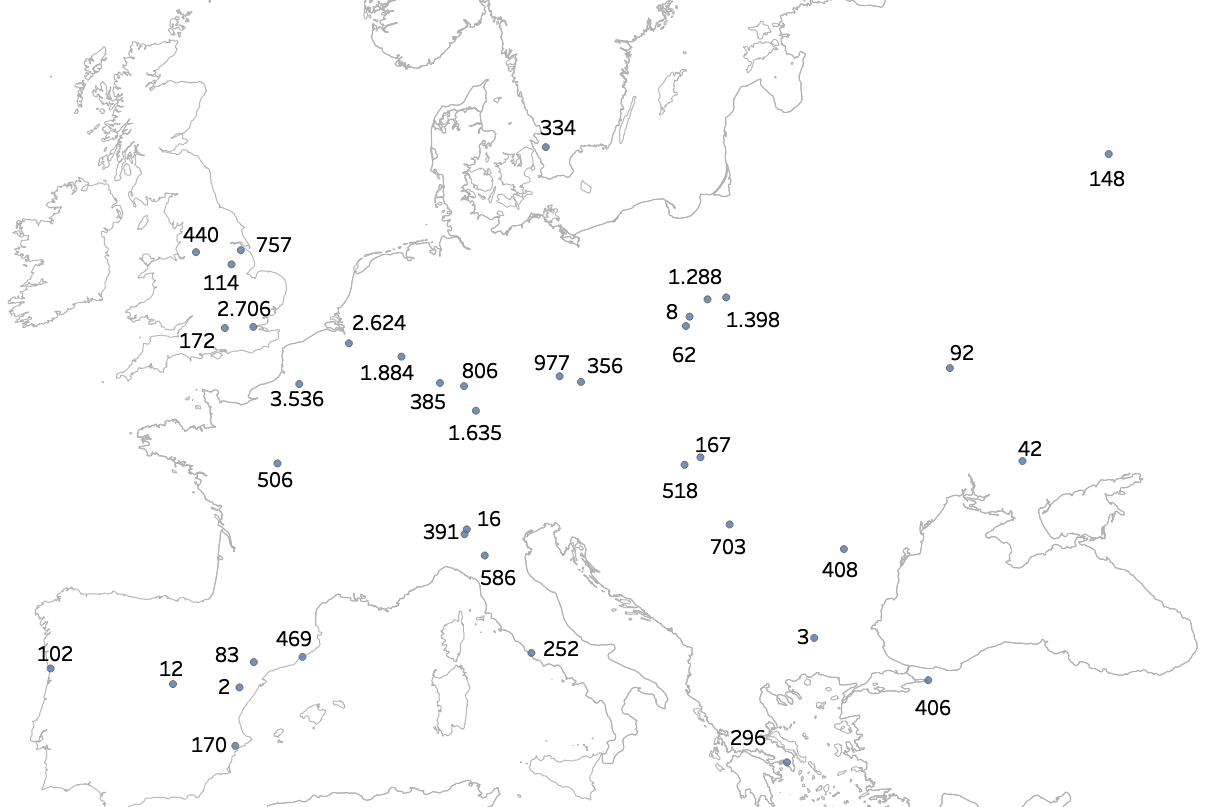}}\fbox{\includegraphics[width=0.48\linewidth]{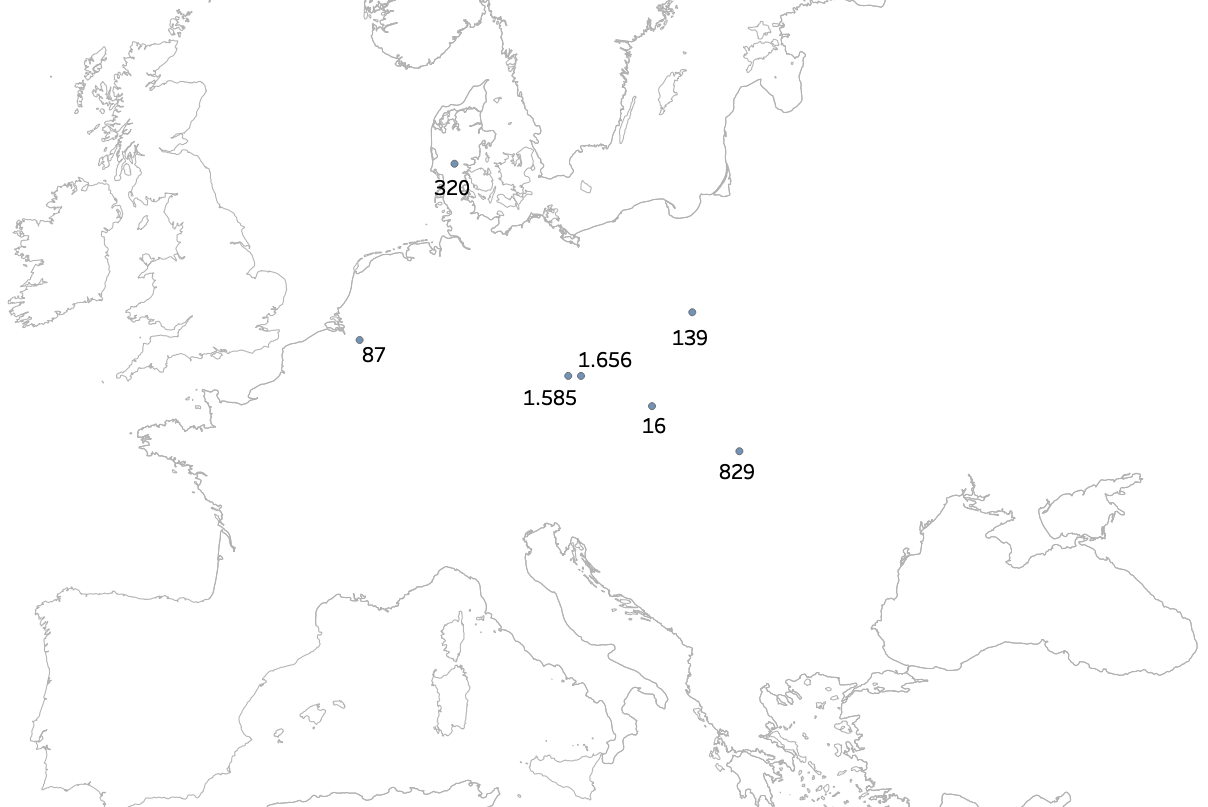}}
\caption{Number of Requests per Location for Customer~1 (left) and 2 (right).\label{fig:ships}}
\end{figure}

\begin{figure}[t]
\centering
\fbox{\includegraphics[width=0.48\linewidth]{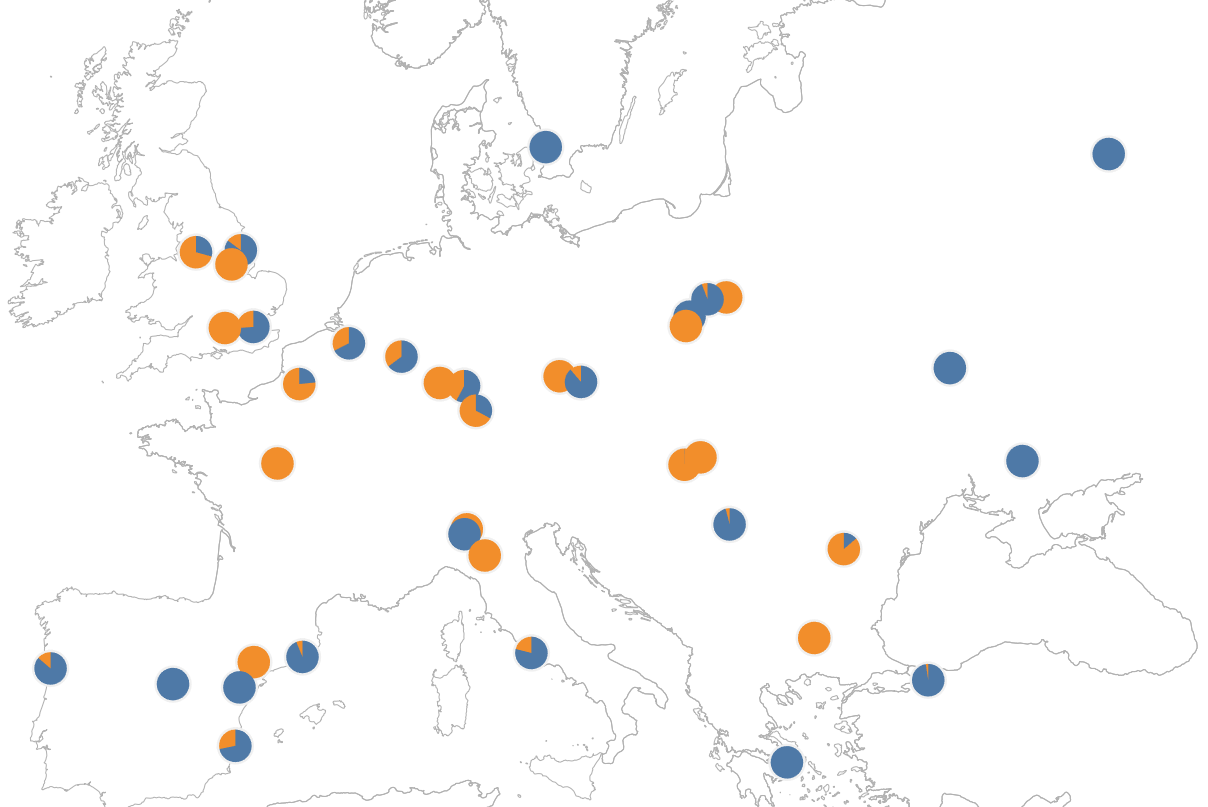}}\fbox{\includegraphics[width=0.48\linewidth]{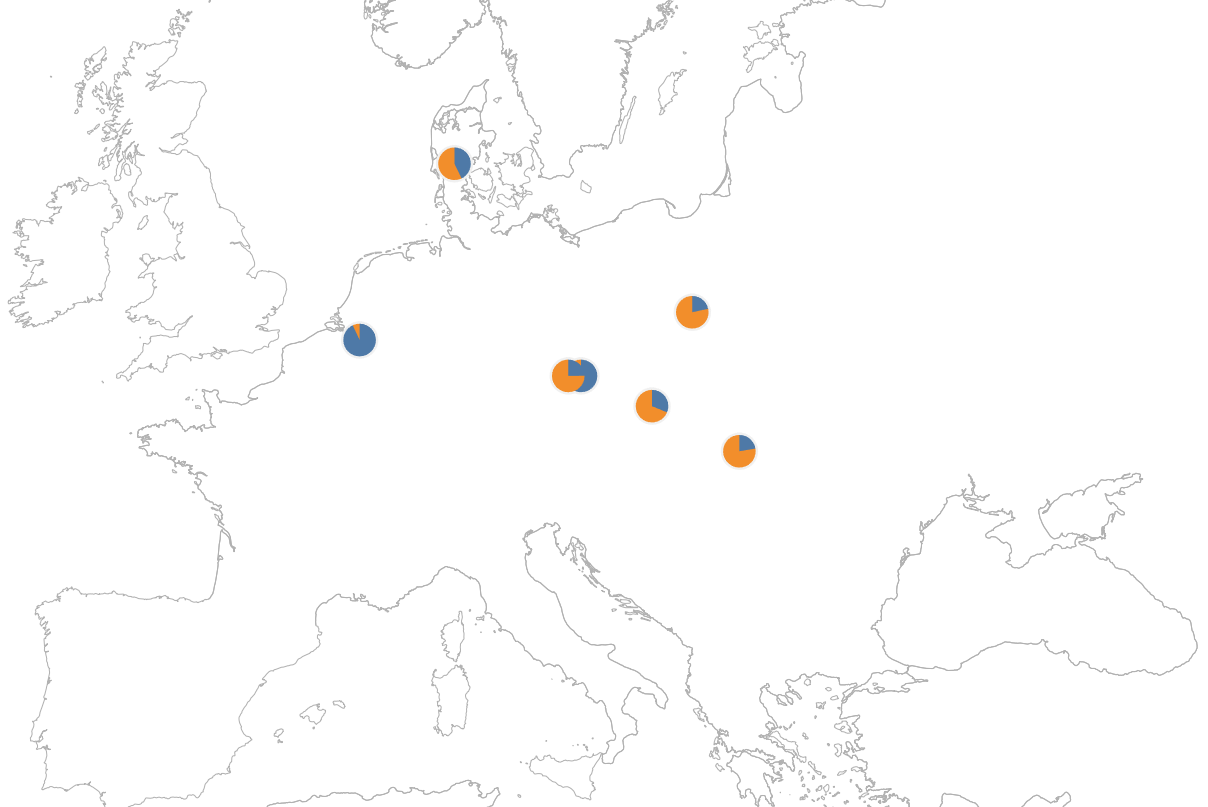}}
\caption{In- and Out-Degree per Location for Customer~1 (left) and 2 (right).\label{fig:mapdegrees}
The pie-chart per location shows the relative in-degree in blue and the out-degree in orange.}
\end{figure}

\paragraph{Request Data}
Two data sets obtained from DB Schenker \emph{customers} are used.
We consider two instances for Customer~1 and three instances for Customer~2, each consisting of monthly request data, representing low-, average- and high-demand months.  
In contrast to Customer~2, the demand for Customer~1 is relatively stable over time. Hence the difference between high- and low-volume months is rather insignificant.
Table~\ref{tab:instances} provides the total number of requests, along with the corresponding sum of the direct origin-destination distances over all requests for all instances.
Figure~\ref{fig:degrees} visualizes the number of in- and out-going requests per location.
Note that most locations have a significant imbalance between in- and out-going requests.
Additionally, the maps in Figures~\ref{fig:ships} and~\ref{fig:mapdegrees} visualize the number of requests as well as the in- and out-degree.

\begin{table}
\centering\footnotesize
\begin{tabular}{c|rr|rr|rr}%
\hline%
&\multicolumn{2}{c|}{High Volume Month}&\multicolumn{2}{c|}{Average\ Volume Month}&\multicolumn{2}{c}{Low Volume Month}\\%
Customer  &\#requests & kilometer & \#requests & kilometer & \#requests & kilometer\\%
\hline%
1&12,427&10,215,308&---&---&11,432&9,197,733\\%
2&2,316&1,094,550&2,257&1,145,214&1,162&678,370\\%
\hline%
\end{tabular}
\caption{Instances}
\label{tab:instances}
\end{table}

Figure \ref{fig:mapdegrees} indicates that for Customer 1, quite some locations are imbalanced, \ie, has either (almost) only incoming, or (almost) only outgoing requests.
In contrast, the demand from Customer 2 looks more balanced, except for the western location.

Distances between physical locations were computed using the Open Source Routing Machine \cite{luxen-vetter-2011}.
    
Overall, it can be observed that the problem instances are large, with up to $12{,}427$ shipments.

As aligned with the FTL operations team, the minimum amount to travel by an own vehicle was set to $\minkm = 8{,}000$\,km and $\minkm = 5{,}000$\,km, for Customer 1 and 2, respectively.

Note that the ultimate goal is to generate insights regarding the \ac{FCT} fleet size for each customer, as operational plans are out of scope in this paper.
Hence, solving multiple scenarios depending on different monthly volumes helps the sales department develop cost-efficient and reliable business offers.

\paragraph{Cost Data}

DB Schenker's FTL operations team provided us with \ac{SM} rates to be used for outsourced requests corresponding to the given datasets.
Since it was impossible to obtain all rates (for all points in time and all origin-destination pairs), we had to make assumptions regarding missing rates.
In close collaboration with our business partners, we managed to develop realistic assumptions for static \ac{SM} rates per kilometer.

\begin{figure}
\centering
\includegraphics[trim=1.2in 7.2in 0.85in 1.05in,clip,width=\linewidth]{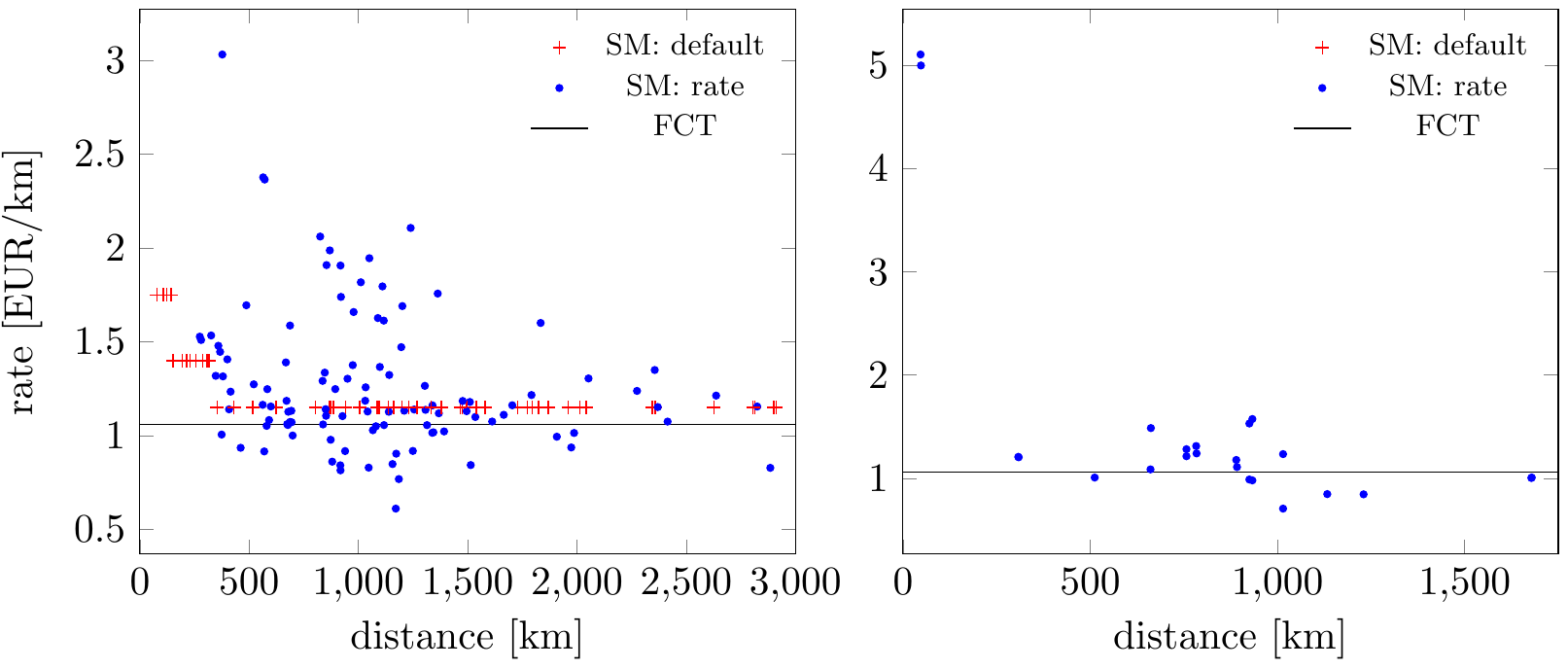}
\caption{Real \ac{SM} rates for out-sourcing vs.\ approximated rates for out-sourcing vs.\ costs per km for vehicles for Customer~1 (left) and 2 (right).\label{fig:rateperkm}
The individual points show the approximated \emph{\ac{SM} default} rates (red crosses) as well as the real \emph{\ac{SM} rate}s for outsourcing (blue dots). The static cost of 1.06\,EUR/km for vehicles is shown as a black reference line.}
\end{figure}

In particular, we distinguish three different specific rate levels, depending on the distance between origin and destination, as shown in Table~\ref{tab:smRates_lit}.
Regarding the cost per driven kilometer of the \ac{FCT} vehicles, our business partners provided a cost benchmark of 1.06\,EUR/km (both empty or loaded).
Figure~\ref{fig:rateperkm} displays the relation between specific rates in EUR/km and the distance from the origin to the destination.
In particular, \emph{\ac{SM}: default} rates indicate the distance-specific approximation described above, whereas \emph{\ac{SM}: rate} indicates the real rates for out-sourcing to the \ac{SM} obtained from the operations team. 

For Customer~1, approximately 35\,\% of the \ac{SM} rates could not be obtained by the operations team and were, thus, approximated using the approach described above.
In contrast to that, all \ac{SM} rates could be retrieved for Customer~2. 
No clear trend can be observed in Figure~\ref{fig:rateperkm}. 
This is mainly because specific rate levels differ strongly between trade lanes, and that they are not symmetric per trade lane: a transport from Eastern Europe to Western Europe, with a distance of, \eg, 750\,km would cost less than the corresponding trip in the opposite direction.
Also, rate levels are dynamic in time.
Certain months of the year experience high demand for transportation services, hence the costs get higher.
Similarly, for low-volume months the transportation costs are expected to be cheaper.

Note that these costs reflect the market at the time of writing this article.

\subsubsection{Out-Sourcing Everything vs.\ Nothing}
\label{subsection:outsourcingeverythingnothing}

\noindent
Here, we present the computational results of enforcing all requests to be outsourced and all requests to be served by vehicles.
Computing the cost of outsourcing is easy. We have to sum up all outsourcing costs of all requests.
For computing the cost of serving all requests by vehicles, we consider the cost of unplanned requests within the \ac{ALNS} as a huge number.

Tables~\ref{tab:computationalResults} and~\ref{tab:businessResults} indicate the results considering complete sets of requests for each customer and month.
In particular, Table~\ref{tab:computationalResults} shows the number of vehicles needed to perform the service.
Additional columns present the empty travel distance between loads relative to the loaded travel distance, \ie, between the previous delivery and the next pickup, KPIs related to the distance traveled per vehicle, and, finally, the computational time needed.

\begin{table}
\centering\footnotesize
\begin{tabular}{ccr|r|rrrr|r}%
\hline%
Cus{-}&&&&Additional&\multicolumn{3}{c|}{Dist. p. Vehicle}&CPU\\%
tomer&Month&\#requests&\#vehicles&Empty Dist.&Max&Avg&Min&{[}sec{]}\\%
\hline%
1 & High & 12,427 & 1,172 & +31.3\,\% & 18.8 & 11.7 & 8.0 &4,386\\%
1 & Low  & 11,432 & 1,095 & +35.1\,\% & 18.8 & 11.9 & 8.0 &3,917\\%
\hline%
2 & High & 2,316 & 153 & +35.6\,\% & 13.2 &  9.7 & 5.4 & 525\\%
2 & Avg  & 2,257 & 153 & +41.1\,\% & 13.8 & 10.6 & 5.4 & 727\\%
2 & Low  & 1,162 &  95 & +33.8\,\% & 12.5 &  9.6 & 5.4 & 345\\%
\hline%
\end{tabular}
\caption{Computational results for a solution out-sourcing nothing.\label{tab:computationalResults} Distances per vehicle are given in $1{,}000$\,km. The additional empty distances are given relatively to the loaded distances in Table~\ref{tab:instances}.}
\end{table}

\begin{table}
\centering\footnotesize
\begin{tabular}{@{}c@{~}c||r@{q\quad}r@{\quad}r@{\quad}r@{\quad}r@{\quad}r|r||r|r@{}}%
\multicolumn{2}{l||}{out-sourcing} & \multicolumn{7}{c||}{everything} & \multicolumn{2}{c}{nothing}\\\hline\hline
      &              & \multicolumn{6}{c|}{requests grouped per orig.-dest.~dist.~$\delta$} &&\\
Cus-    &       & \multicolumn{2}{c}{$\delta{<}150$} & \multicolumn{2}{c}{$150{\leq} \delta{<}350$} & \multicolumn{2}{c|}{$350{\leq} \delta$} &&\\
tomer & Month & \#req. & $\Sigma$km & \#req. & $\Sigma$km & \#req. & $\Sigma$km & cost & km & cost\\
\hline
1&High&458&46.9&2,651&726.3&9,318&9,442.1&\bf 12.43&13,417&14.22\\%
1&Low &435&44.3&2,650&747.4&8,347&8,406.0&\bf 11.40&12,431&13.18\\%
\hline%
2&High&1,099&52.7&16&4.9&1,201&1,036.9&\bf 1.50&1,485&1.57\\%
2&Avg &  958&45.9&24&7.4&1,275&1,091.9&\bf 1.50&1,616&1.71\\%
2&Low &  439&21.0&17&5.3&  706&  652.1&\bf   0.85&  908&  0.96\\%
\hline
\end{tabular}
\caption{Costs for solutions out-sourcing everything and out-sourcing nothing.\label{tab:businessResults} Distances are given in $1{,}000$\,km. Costs are given in $1$ million EUR. The highlighted cost is better. \#req.\ denotes the number of requests.}
\end{table}

It can easily be seen that significant additional empty travel would be required in a scenario where nothing is outsourced.
One of the main reasons for this is, of course, the imbalanced demand structure that was displayed in Figure~\ref{fig:mapdegrees}.
In particular, several locations in both data sets serve either as origins (sources) or destinations (sinks). This naturally leads to empty runs, \eg, when a vehicle delivers at a pure sink location and then is forced to travel empty to its following source location.

To assess the merits of a scenario without outsourcing, we compare its overall cost to that of a  scenario where everything is outsourced, as summarized in Table~\ref{tab:businessResults}.
In all five instances, outsourcing nothing is less competitive than outsourcing everything, with an overall cost disadvantage of up to 16\,\%.
Again, this result does not come as a real surprise since an outsourcing approach is, at least in theory, fully flexible and of unlimited capacity and, thus, more efficient in serving the rather imbalanced and erratic demand patterns inherent to our customer instances.

\subsubsection{Cost-Optimal Mix}

After computing the costs for outsourcing no request and outsourcing all requests for comparison, we compute the cost-optimal assignment of requests to outsourcing options or an own fleet of vehicles.
We call a scenario where both options are allowed a \emph{mixed scenario}.

Table~\ref{tab:subset} indicates results for the most cost-efficient mixed scenarios found and compares them to outsourcing nothing and outsourcing everything. In all five instances, the mixed scenario is more competitive than any other two options, although the relative gain is more marked for Customer~2 than Customer~1.

\begin{table}
\centering\footnotesize
\begin{tabular}{@{}l|r@{~~}r|r@{~~}r|r@{~~}r@{~~}r@{~~}r@{~}|@{~}r@{~~}r|r@{}}
\hline
Cus- & \multicolumn{2}{c|}{out-sourcing} & \multicolumn{2}{c|}{out-sourcing} & \multicolumn{7}{c}{mixed scenario}\\
tomer & \multicolumn{2}{c|}{nothing} & \multicolumn{2}{c|}{everything} & \multicolumn{4}{c@{~}|@{~}}{own vehicles} & \multicolumn{2}{c|}{out-sourced} & $\Sigma$\\
\&  & \multicolumn{2}{c|}{} & \multicolumn{2}{c|}{} & \% & km & km & &&& \\
Month & km & cost & km & cost & req. & loaded & empty & cost & km & cost & cost\\
\hline
1\,High & 13,417 & 14.22 & 10,215 & 12.43 &     6.9 &   294.4 &   39.3 &  0.35 & 9,921 & 11.83 & \bf 12.18 \\
1\,Low  & 12,431 & 13.18 &  9,198 & 11.40 &     6.2 &   264.8 &   24.2 &  0.31 & 8,933 & 10.87 & \bf 11.18 \\\hline
2\,High &  1,485 &  1.57 &  1,095 &  1.50 &    41.1 &    56.2 &   35.7 &  0.10 & 1,038 &  1.26 & \bf 1.36 \\
2\,Avg  &  1,616 &  1.71 &  1,145 &  1.50 &    35.4 &    50.2 &   26.0 &  0.08 & 1,095 &  1.30 & \bf 1.38 \\
2\,Low  &    908 &  0.96 &    678 &  0.85 &    34.5 &    26.3 &   15.0 &  0.04 &   652 &  0.74 & \bf 0.79 \\\hline
\end{tabular}
\caption{Costs for the mixed scenarios. Distances are given in $1{,}000$\,km. \%req.\ denotes the percentage of requests. Costs are given in $1$ million EUR. The highlighted cost is better.\label{tab:subset}}
\end{table}

The demand structure of Customer~1 (as discussed in Section~\ref{subsection:datadescription}) shows an extreme imbalance of in- and out-going requests for many locations.
The geographic distribution of the locations induces many long-distance connections, which lead to a high amount of empty travel.
To be financially beneficial, a spot-market rate has to be twice as high as the rate for operating an own vehicle.
This is rarely the case, as shown before.
As a result of this effect, we see that only a minimal amount (up to 7\,\%) of the requests is served by own vehicles, whereas the vast majority is outsourced.
For a high-volume month, the mixed scenario is 2\,\% cheaper (12.18 million instead of 12.43 million EUR) than outsourcing everything to the \ac{SM}.
In a low-volume month, savings are even less.

As the demand structure for Customer~2 is more balanced, significantly less empty travel is induced.
This results in serving up to 41\,\% of the requests using own vehicles.
For a high/average/low-volume month, the mixed scenario is 9\,\%/7\,\%/6\,\% cheaper than out-sourcing everything to the spot-market, respectively.

\section{Conclusions}

\noindent
In this paper, we have presented a meta-heuristic to effectively support business development units at DB Schenker in designing competitive offers for complex full-load solutions. 
In particular, in order to calculate the total number of vehicles needed, the distance traveled, and the full-load requests that are served by an own fleet of vehicles instead of out sourcing them to \ac{SM}. We modeled the problem as a variant of the \ac{PDPTW} with driving regulations, and tackled it using a tailored \ac{ALNS}. 

We compared three scenarios, namely out-sourcing nothing, out-sourcing everything to the spot-market, and a mix of both, and quantified the corresponding costs.
For evaluation purposes, we used transformed VRP instances widely used in the scientific literature (\ie, Gehring \& Homberger) and real-life instances from two potential customers of DB Schenker, containing monthly demand data with up to 12,427 requests.
For these instances, out-sourcing nothing is outperformed by out-sourcing everything due to the underlying demand structure. 
However, a mixed setup may yield benefits of up to 9\,\% compared to out-sourcing everything.

To actually realize the potential benefits presented in this paper, the developed approach would need to be complemented by an operational decision support system, which can help control towers continuously plan, execute, and re-plan scenarios in case of demand fluctuations, driving time deviations or other unexpected events.
The operational system needs access to real-time traffic data and own vehicles should be equipped with tracking devices that allow timely detection of potential deviations or disruptions.

Also from an algorithmic point of view, our approach could be extended in a number of ways.
As aforementioned, when using the presented algorithm as operational decision support, additional aspects/parameters need to be added to allow for considering dynamic aspects of the problem as well as problem heterogeneity, \eg, different cost assumptions per lane or geography when operating own vehicles. 
We have also indicated above that we took a rather pragmatic approach towards driving time constraints which currently doesn't reflect the full complexity of various national regulatory regimes.

In addition, when using the approach for tender calculations, demand uncertainty should be incorporated in the optimization procedure. 
This implies that routing solutions would need to be more conservative, hence more costly to some extent.
However, incorporating the uncertainty explicitly into the model would increase the robustness of the solutions, and thus lead to less probability of unexpected costs due to demand fluctuations.

\bibliography{mybibfile}

\begin{thebibliography}{10}
\expandafter\ifx\csname url\endcsname\relax
  \def\url#1{\texttt{#1}}\fi
\expandafter\ifx\csname urlprefix\endcsname\relax\def\urlprefix{URL }\fi
\expandafter\ifx\csname href\endcsname\relax
  \def\href#1#2{#2} \def\path#1{#1}\fi

\bibitem{wieberneit}
N.~Wieberneit, Service network design for freight transportation: A review, OR
  Spectrum 30 (2008) 77--112.
\newblock \href {https://doi.org/10.1007/s00291-007-0079-2}
  {\path{doi:10.1007/s00291-007-0079-2}}.

\bibitem{sav_sol1995}
M.~Sol, M.~W.~P. Savelsbergh, The general pickup and delivery problem,
  Transportation Science 29 (1995) 17--29.

\bibitem{SOARES2019174}
R.~Soares, A.~F. Marques, P.~Amorim, J.~Rasinm{\"{a}}ki, Multiple vehicle
  synchronisation in a full truck-load pickup and delivery problem: {A}
  case-study in the biomass supply chain, European Journal of Operational
  Research 277~(1) (2019) 174--194.
\newblock \href {https://doi.org/10.1016/j.ejor.2019.02.025}
  {\path{doi:10.1016/j.ejor.2019.02.025}}.

\bibitem{XUE2021500}
N.~Xue, R.~Bai, R.~Qu, U.~Aickelin, A hybrid pricing and cutting approach for
  the multi-shift full truckload vehicle routing problem, European Journal of
  Operational Research 292~(2) (2021) 500--514.
\newblock \href {https://doi.org/10.1016/j.ejor.2020.10.037}
  {\path{doi:10.1016/j.ejor.2020.10.037}}.

\bibitem{asvinGoel2009}
A.~Goel, {Vehicle Scheduling and Routing with Drivers' Working Hours},
  Transportation Science 43~(1) (2009) 17--26.
\newblock \href {https://doi.org/10.1287/trsc.1070.0226}
  {\path{doi:10.1287/trsc.1070.0226}}.

\bibitem{Goel:2010:TDS:1886581.1886582}
A.~Goel, {Truck Driver Scheduling in the European Union}, Transportation
  Science 44~(4) (2010) 429--441.
\newblock \href {https://doi.org/10.1287/trsc.1100.0330}
  {\path{doi:10.1287/trsc.1100.0330}}.

\bibitem{GOEL2018144}
A.~Goel, {Legal aspects in road transport optimization in Europe},
  Transportation Research Part E: Logistics and Transportation Review 114
  (2018) 144--162.
\newblock \href {https://doi.org/10.1016/j.tre.2018.02.011}
  {\path{doi:10.1016/j.tre.2018.02.011}}.

\bibitem{GOEL2013a}
A.~Goel, T.~Vidal, Hours of service regulations in road freight transport: An
  optimization-based international assessment, Transportation Science 48 (2013)
  313--463.
\newblock \href {https://doi.org/10.1287/trsc.2013.0477}
  {\path{doi:10.1287/trsc.2013.0477}}.

\bibitem{GOEL2016a}
A.~Goel, S.~Irnich, An exact method for vehicle routing and truck driver
  scheduling problems, Transportation Science 51 (2016) 395--789.
\newblock \href {https://doi.org/10.1287/trsc.2016.0678}
  {\path{doi:10.1287/trsc.2016.0678}}.

\bibitem{palgrave.jors.2602143}
P.~P. Repoussis, C.~D. Tarantilis, G.~Ioannou, The open vehicle routing problem
  with time windows, Journal of the Operational Research Society 58~(3) (2007)
  355--367.
\newblock \href {https://doi.org/10.1057/palgrave.jors.2602143}
  {\path{doi:10.1057/palgrave.jors.2602143}}.

\bibitem{doi:10.1080/00207543.2019.1572929}
R.~Lahyani, A.-L. Gouguenheim, L.~C. Coelho, A hybrid adaptive large
  neighbourhood search for multi-depot open vehicle routing problems,
  International Journal of Production Research 57~(22) (2019) 6963--6976.
\newblock \href {https://doi.org/10.1080/00207543.2019.1572929}
  {\path{doi:10.1080/00207543.2019.1572929}}.

\bibitem{RePEc:eee:ejores:v:231:y:2013:i:1:p:1-21}
T.~Vidal, T.~G. Crainic, M.~Gendreau, C.~Prins, Heuristics for multi-attribute
  vehicle routing problems: a survey and synthesis, European Journal of
  Operational Research 231~(1) (2013) 1--21.
\newblock \href {https://doi.org/10.1016/j.ejor.2013.02.053}
  {\path{doi:10.1016/j.ejor.2013.02.053}}.

\bibitem{Braekers2016TheVR}
K.~Braekers, K.~Ramaekers, I.~V. Nieuwenhuyse, The vehicle routing problem:
  State of the art classification and review, Computers \& Industrial
  Engineering 99 (2016) 300--313.
\newblock \href {https://doi.org/10.1016/j.cie.2015.12.007}
  {\path{doi:10.1016/j.cie.2015.12.007}}.

\bibitem{lit_review2016}
A.~Annouch, K.~Bouyahyaoui, A.~Bellabdaoui, A literature review on the full
  trackload vehicle routing problems, in: A.~E.~H. Alaoui, Y.~Benadada,
  J.~Boukachour (Eds.), 3rd International Conference on Logistics Operations
  Management, {GOL} 2016, Fez, Morocco, May 23-25, 2016, IEEE, 2016, pp. 1--6.
\newblock \href {https://doi.org/10.1109/GOL.2016.7731723}
  {\path{doi:10.1109/GOL.2016.7731723}}.

\bibitem{doi:10.1057/palgrave.jors.2602603}
C.~Archetti, D.~Feillet, A.~Hertz, M.~G. Speranza, The capacitated team
  orienteering and profitable tour problems, Journal of the Operational
  Research Society 60~(6) (2009) 831--842.
\newblock \href {https://doi.org/10.1057/palgrave.jors.2602603}
  {\path{doi:10.1057/palgrave.jors.2602603}}.

\bibitem{orient_pdp}
M.~Gansterer, M.~Küçüktepe, R.~Hartl, The multi-vehicle profitable pickup
  and delivery problem, OR Spectrum 39 (06 2016).
\newblock \href {https://doi.org/10.1007/s00291-016-0454-y}
  {\path{doi:10.1007/s00291-016-0454-y}}.

\bibitem{sun}
P.~Sun, L.~Veelenturf, M.~Hewitt, T.~Van~Woensel, Adaptive large neighborhood
  search for the time-dependent profitable pickup and delivery problem with
  time windows, Transportation Research Part E: Logistics and Transportation
  Review 138 (2020) 101942.
\newblock \href {https://doi.org/10.1016/j.tre.2020.101942}
  {\path{doi:10.1016/j.tre.2020.101942}}.

\bibitem{Kara2011}
I.~Kara, Arc based integer programming formulations for the distance
  constrained vehicle routing problem, in: 3rd IEEE International Symposium on
  Logistics and Industrial Informatics, 2011, pp. 33--38.
\newblock \href {https://doi.org/10.1109/LINDI.2011.6031159}
  {\path{doi:10.1109/LINDI.2011.6031159}}.

\bibitem{GoelGruhn2006}
A.~Goel, V.~Gruhn, Drivers' working hours in vehicle routing and scheduling,
  in: {IEEE} Intelligent Transportation Systems Conference, {ITSC} 2006,
  Toronto, Ontario, Canada, 17-20 September 2006, {IEEE}, 2006, pp. 1280--1285.
\newblock \href {https://doi.org/10.1109/ITSC.2006.1707399}
  {\path{doi:10.1109/ITSC.2006.1707399}}.

\bibitem{Shaw1998}
P.~Shaw, Using constraint programming and local search methods to solve vehicle
  routing problems, in: M.~J. Maher, J.~Puget (Eds.), Principles and Practice
  of Constraint Programming -- CP98, 4th International Conference, Pisa, Italy,
  October 26-30, 1998, Proceedings, Vol. 1520 of Lecture Notes in Computer
  Science, Springer, 1998, pp. 417--431.
\newblock \href {https://doi.org/10.1007/3-540-49481-2_30}
  {\path{doi:10.1007/3-540-49481-2_30}}.

\bibitem{Pisinger2007ALNS}
D.~Pisinger, S.~Ropke, A general heuristic for vehicle routing problems,
  Computers \& Operations Research 34~(8) (2007) 2403--2435.
\newblock \href {https://doi.org/10.1016/j.cor.2005.09.012}
  {\path{doi:10.1016/j.cor.2005.09.012}}.

\bibitem{ropke2006}
S.~Ropke, D.~Pisinger, An adaptive large neighborhood search heuristic for the
  pickup and delivery problem with time windows, Transportation Science 40
  (2006) 455--472.
\newblock \href {https://doi.org/10.1287/trsc.1050.0135}
  {\path{doi:10.1287/trsc.1050.0135}}.

\bibitem{Ghilas:2016:ALN:2927990.2928062}
V.~Ghilas, E.~Demir, T.~Van~Woensel, An adaptive large neighborhood search
  heuristic for the pickup and delivery problem with time windows and scheduled
  lines, Computers \& Operations Research 72 (2016) 12--30.
\newblock \href {https://doi.org/10.1016/j.cor.2016.01.018}
  {\path{doi:10.1016/j.cor.2016.01.018}}.

\bibitem{GRIMAULT20171}
A.~Grimault, N.~Bostel, F.~Lehu{\'{e}}d{\'{e}}, An adaptive large neighborhood
  search for the full truckload pickup and delivery problem with resource
  synchronization, Computers \& Operations Research 88 (2017) 1--14.
\newblock \href {https://doi.org/10.1016/j.cor.2017.06.012}
  {\path{doi:10.1016/j.cor.2017.06.012}}.

\bibitem{potvin1993}
J.-Y. Potvin, J.-M. Rousseau, A parallel route building algorithm for the
  vehicle routing and scheduling problem with time windows, European Journal of
  Operational Research 66~(3) (1993) 331--340.
\newblock \href {https://doi.org/10.1016/0377-2217(93)90221-8}
  {\path{doi:10.1016/0377-2217(93)90221-8}}.

\bibitem{savelsbergh1990}
M.~W.~P. Savelsbergh, An efficient implementation of local search algorithms
  for constrained routing problems, European Journal of Operational Research
  47~(1) (1990) 75--85.
\newblock \href {https://doi.org/10.1016/0377-2217(90)90091-O}
  {\path{doi:10.1016/0377-2217(90)90091-O}}.

\bibitem{campbell2004}
A.~M. Campbell, M.~W.~P. Savelsbergh, Efficient insertion heuristics for
  vehicle routing and scheduling problems, Transportation Science 38 (2004)
  369--378.
\newblock \href {https://doi.org/10.1287/trsc.1030.0046}
  {\path{doi:10.1287/trsc.1030.0046}}.

\bibitem{GehringHomberger1999}
H.~Gehring, J.~Homberger,
  \href{http://www.mit.jyu.fi/eurogen99/papers/homberg.ps}{A parallel hybrid
  evolutionary metaheuristic for the vehicle routing problem with time
  windows}, in: Proceedings of the EUROGEN99: Short Course on Evolutionary
  Algorithms in Engineering and Computer Science, University of
  Jyv\"{a}skyl\"{a}, Finland, 1999.
\newline\urlprefix\url{http://www.mit.jyu.fi/eurogen99/papers/homberg.ps}

\bibitem{luxen-vetter-2011}
D.~Luxen, C.~Vetter, Real-time routing with openstreetmap data, in: Proceedings
  of the 19th ACM SIGSPATIAL International Conference on Advances in Geographic
  Information Systems, GIS '11, ACM, New York, NY, USA, 2011, pp. 513--516.
\newblock \href {https://doi.org/10.1145/2093973.2094062}
  {\path{doi:10.1145/2093973.2094062}}.

\end{thebibliography}

\end{document}